\UseRawInputEncoding
\newtheorem{remark}{Remark}

\newtheorem{theorem}{Theorem}
\newtheorem{lemma}{Lemma}

\documentclass[10pt,journal,compsoc]{IEEEtran}
\newcommand{\bea}{\begin{eqnarray}}
\newcommand{\eea}{\end{eqnarray}}

\usepackage[dvips]{graphicx}

\usepackage{amsmath,amssymb,amsfonts}
\usepackage{cite}
\usepackage{graphicx}
\usepackage{subfigure}
\usepackage{color}
\renewcommand{\baselinestretch}{1}
\usepackage{geometry}
\begin{document}

\title{An Elegant Inequality }

\author{Yiguang Liu
\IEEEcompsocitemizethanks{\IEEEcompsocthanksitem Yiguang Liu is with College of Computer Science and with School of Aeronautics and Astronautics, Sichuan University, Chengdu, Sichuan
Province, China, 610065. E-mail: liuyg@scu.edu.cn
}}
\maketitle
\begin{abstract}
A new inequality, $(x)^{p}+(1-x)^{\frac{1}{p}}\leq1$ for $p \geq 1$ and $\frac{1}{2} \geq x \geq 0$ is found and proved.
The inequality looks elegant as it integrates two number pairs ($x$ and $1-x$, $p$ and $\frac{1}{p}$) whose summation and product are one. Its right hand side, $1$, is the strict upper bound of the left hand side. The equality cannot be categorized into any known type of inequalities such as H\"{o}lder, Minkowski etc. In proving it, transcendental equations have been met with, so some novel techniques have been built to get over the difficulty.
\end{abstract}


\section{Introduction}
\IEEEPARstart{I}{nequalities} are critical for studying performance of algorithms, models and others. For example, inequalities involving eigenvalues and singular values of products and sums of matrices are useful to analyze multi-variable systems \cite{571}; Sobolev Type Inequalities have a strong influence on development of the theory of partial differential equations, analysis, mathematical physics, differential geometry, and other fields \cite{552}. So, many books about inequalities have been published. For instance, \cite{1082} consists of more than 5000 inequalities and 50 methods for proving inequalities.

Recently, inequalities still keep rapid developing. Correlation inequalities have been particularly paid attention to for studying stochastic models for large interacting systems \cite{489}.  Bernoullis inequality is often used as the crucial step in the proof of other inequalities, and various generalizations of it were given in \cite{1081}. As is well known, the H\"{o}lder's inequality has important applications in many areas of pure and applied mathematics, and a new sharpened and generalized version of H\"{o}lder's inequality has been given in \cite{1083}. It has also been refined and extended to the case of multiple sequences \cite{1084}. New generalizations of Acz\'{e}l's inequality and Popoviciu's inequality were given \cite{1085} , which has broad applications to mathematical analysis. In \cite{1087}, a generalization of the Ky Fan discrete multivariate inequality was obtained and some of its applications were also given.

Though there are so many classical inequalities, such as fundamental inequalities, combinatorial inequalities, variational inequalities, determinant and matrix inequalities, sequence and series inequalities, differential inequalities, etc, we don't meet with an algebraic inequality which consists of two numbers which are reciprocal, and of another two numbers whose summation is $1$. In this paper, we put forward a formally elegant inequality having this property.

\section{Mathematical Apparatus} \label{sec.ma}
\subsection{Primaries}
In calculus, a derivative test is used to locate the critical points of a function and determine whether such points are a local maximum or a local minimum. The usefulness of derivatives to find extrema is proved mathematically by Fermat's theorem of stationary points. To make this paper self-contained, some criteria borrowed from textbooks such as \cite{532} are listed as follows.
\begin{lemma} \label{1.lemma} Let $f:(a,b)\rightarrow \mathbf{R}$ be a function and suppose that $x_{0} \in (a ,b )$ is a point where $f$ has a local extremum. If $f$ is differentiable at $x_{0}$, then $\dot{f}(x_{0})=0$. That is,  if the derivative of a function at any point is not zero, then there is not a local extremum at that point.
\end{lemma}
\begin{lemma} \label{2.lemma} If $x_{0}$ is an extremum of $f$, then one of the following is true: 1) $x_{0}$ is at $a$ or $b$; 2) $f$ is not differentiable at $x_{0}$; 3) $x_{0}$ is a stationary point of $f$.
\end{lemma}
\begin{lemma} \label{3.lemma} If the function $f$ is twice-differentiable at stationary point $x_{0}$, then: 1) if $\ddot{f}(x_{0})<0$, then $f$ attains a local maximum at $x_{0}$; 2) if $\ddot{f}(x_{0})>0$, then $f$ attains a local minimum at $x_{0}$; 3) if $\ddot{f}(x_{0})=0$, then whether $f(x_{0})$ is a local extreme is inconclusive.
\end{lemma}

\subsection{Main Contributions}
Using derivative test, we have proved the proposed inequality in this section. First, we prove the following theorem.
\begin{theorem} \label{1.theo}  For $p\geq1$, $(\frac{1}{2})^{p}+(\frac{1}{2})^{\frac{1}{p}}\leq1$.
\end{theorem}
 {\it {\bf Proof}}: Let
\bea \label{5.eq.lemma1}
f(p)=(\frac{1}{2})^{p}+(\frac{1}{2})^{\frac{1}{p}}, p\geq1.
\eea From \eqref{5.eq.lemma1} we can see $f(p)$ is always derivative, so there are
\bea \label{6.eq.lemma1}
\dot{f}(p)=-(\frac{1}{2})^{p}\log(2)+\frac{1}{p^2}(\frac{1}{2})^{\frac{1}{p}}\log(2),
\eea and
\bea \label{7.eq.lemma1}
\ddot{f}(p)=(\frac{1}{2})^{p}\log^{2}(2)-2\frac{1}{p^3}(\frac{1}{2})^{\frac{1}{p}}\log(2)+\frac{1}{p^4}(\frac{1}{2})^{\frac{1}{p}}\log^{2}(2).
\eea
At stationary points, $\dot{f}(p)=0$, from \eqref{6.eq.lemma1} there is
\bea \label{8.eq.lemma1}
(\frac{1}{2})^{p}=\frac{1}{p^2}(\frac{1}{2})^{\frac{1}{p}}\equiv C>0.
\eea
To determine whether $f$ is locally maximal or locally minimal, we need to consider the twice derivative of $f$. So substituting \eqref{8.eq.lemma1} into \eqref{7.eq.lemma1}, we have
\bea
\ddot{f}(p)=C(\log^{2}(2)-2\frac{1}{p}\log(2)+\frac{1}{p^2}\log^{2}(2))
\\ \nonumber
=C(\frac{1}{p}-t_{1})(\frac{1}{p}-t_{2})\hspace{2.7cm}
\eea
where $t_{1}\equiv \frac{1-\sqrt{1-\log^2(2)}}{\log(2)}\approx 0.4028<1$, $t_{2}\equiv\frac{1+\sqrt{1-\log^2(2)}}{\log(2)}\approx2.4826$. In terms of the properties of quadratic functions, and $p\geq 1$, we have
\bea \nonumber
\ddot{f}(p)\left\{
              \begin{array}{l}
                \leq 0 \hbox{ when } \frac{1}{p}\in[t_{1},1]\\
                > 0 \hbox{ when } \frac{1}{p}\in(0,t_{1})
              \end{array}
            \right..
\eea
That is, when $p \in [1,t_{1}^{-1}]\approx[1,2.4826]$, there is $\ddot{f}(p)\leq 0$, which means $f_{p}$ take local maximal values at the stationary points $\dot{f}(p)=0$; when $p >t_{1}^{-1}$, $\ddot{f}(p)>0$ means $f_{p}$ take local minimal values at $\dot{f}(p)=0$. Because we only consider the maximal value of function $f(p)$, so we only consider the stationary points $\dot{f}(p)=0$ for $p \in [1,t_{1}^{-1}]$. In the following, we will prove that there is no stationary point for $p \in [1,t_{1}^{-1}]$ except for $p=1$.

The equation \eqref{8.eq.lemma1} amounts to $2\log(p)=(p-\frac{1}{p})\log(2)$, to discuss whether this equation has solutions within $(1,t_{1}^{-1}]$, we let
\bea \label{9.eq.lemma1}
h(p)=2\log(p)-(p-\frac{1}{p})\log(2).
\eea So we have
\bea \label{10.eq.lemma1}
\dot{h}(p)=\frac{2}{p}-(1+\frac{1}{p^2})\log(2)=-(\frac{1}{p}-t_{1})(\frac{1}{p}-t_{2}).
\eea Combining \eqref{9.eq.lemma1} and \eqref{10.eq.lemma1} tells that: 1) $h(p)=0$ and $\dot{h}(p)=0$ when $p=1$; 2) $h(p)\neq0$ for $p\in(1,t_{1}^{-1}]$, because $h(p)$ keeps decreasing from $h(1)=0$ as $\dot{h}(p)$ keeps negative when $p$ changes from $1$ to $t_{1}^{-1}$. So $\dot{f}(p)=0$ does not have solutions within interval $[1,t_{1}^{-1}]$ except for $p=1$. In addition, from \eqref{7.eq.lemma1} we have $\ddot{f}(1)=\log^2(2)-\log(2)<0$. Based on Lemma \ref{1.lemma}, \ref{2.lemma} and \ref{3.lemma}, it is concluded that $f(p)$ takes the maximal value at $p=1$, so we have
\bea \label{11.eq.lemma1}
f(p)=(\frac{1}{2})^{p}+(\frac{1}{2})^{\frac{1}{p}}\leq f(1)=1.
\eea
This theorem is proved. $\blacksquare$
\begin{theorem} \label{2.theo}  For $p\geq1$ and $0\leq x \leq \frac{1}{2}$, $(x)^{p}+(1-x)^{\frac{1}{p}}\leq1$.
\end{theorem}
{\it {\bf Proof}}: Let $f(x)=x^{p}+(1-x)^{\frac{1}{p}}$. We have
\bea \label{eq.f1deriv}
\dot{f}(x)=px^{p-1}-\frac{1}{p}(1-x)^{\frac{1}{p}-1}
\eea and
\bea \label{eq.f2deriv}
\ddot{f}(x)=p(p-1)x^{p-2}+\frac{1}{p}(\frac{1}{p}-1)(1-x)^{\frac{1}{p}-2}.
\eea From derivative test, we know that the maximal value of $f(x)$ possibly lie at the end points of $[0,\frac{1}{2}]$ or at the stationary points within this interval. So we have two cases: 1) if $\dot{f}(x)\neq 0$ for all $x\in [0,\frac{1}{2}]$, then the maximal value of $f(x)$ only lies at the end points, so
\bea \nonumber
f(x)=x^{p}+(1-x)^{\frac{1}{p}}\leq \max(f(0),f(\frac{1}{2}))
\\ \label{1.eq.lemma1}
=\max(1,(\frac{1}{2})^{p}+(\frac{1}{2})^{\frac{1}{p}}).
\eea 2) if there are stationary points within $[0,\frac{1}{2}]$, there is $\dot{f}(x)=0$, so based on \eqref{eq.f1deriv} we have
\bea \label{2.eq.lemma1}
px^{p-1}=\frac{1}{p}(1-x)^{\frac{1}{p}-1}\equiv C >0.
\eea At the stationary points, substituting \eqref{2.eq.lemma1} into \eqref{eq.f2deriv} we have
\bea \label{3.eq.lemma1}
\ddot{f}(x)=C(p-1)\left(\frac{1}{x}-\frac{1}{p(1-x)}\right).
\eea Because $x\in [0,\frac{1}{2}]$ and $p\geq 1$, based on \eqref{3.eq.lemma1} there are
\bea \label{4.eq.lemma1}
\ddot{f}(x)\left\{
              \begin{array}{ll}
                =0, \hbox{when } p=1,x=\frac{1}{2}; \\
                >0, \hbox{otherwise.}
              \end{array}
            \right.
\eea
In terms of Lemma \ref{3.lemma}, \eqref{4.eq.lemma1} tells that: 1) when $\ddot{f}(x)=0$, the evaluations, $p=1$ and $x=\frac{1}{2}$, make $f(x)=1$. 2) When $\ddot{f}(x)>0$, the stationary points resulted from \eqref{2.eq.lemma1} corresponds to the local minimal extrema. That is to say, in this case the maximal values of $f(x)$ only lie at the endpoints of $[0,\frac{1}{2}]$, because $f(x)$ is continuous and its maximal values cannot take at the stationary points where the twice derivative of $f(x)$ is positive. Combining the two cases tells that the maximal value of $f(x)$ is at the endpoints of $[0,\frac{1}{2}]$ or is $1$, which has also been described by \eqref{1.eq.lemma1}.
Based on Theorem \ref{1.theo}, \eqref{1.eq.lemma1} means this theorem.$\blacksquare$

Theorem \ref{2.theo} has many variants, such as the followings.
\begin{remark} \label{1.rmk}  $(\frac{1}{n})^{p}+(1-\frac{1}{n})^{\frac{1}{p}}\leq1$, for $p\geq1$ and $n\in \mathbf{Z}$ and $n\geq 2$.
\end{remark}
\begin{remark} \label{2.rmk}  $(\frac{1-x}{2})^{p}+(\frac{1+x}{2})^{\frac{1}{p}}\leq1$, for $p\geq1$ and $1\geq x\geq 0$.
\end{remark}
\begin{remark} \label{3.rmk}  $\sin^{2p}(\alpha)+\cos^{\frac{2}{p}}(\alpha)\leq1$, for $p\geq1$ and $\alpha\in[0,\frac{\pi}{4}]$.
\end{remark}
\subsection{Comparisons and Discussions}
In \cite{1083}, a generalized H\"{o}lder inequality is introduced as follows
\bea \nonumber
(a_{11}^{\lambda_{1}}a_{12}^{\lambda_{2}}+a_{21}^{\lambda_{1}}a_{22}^{\lambda_{2}})2^{\min(0,\lambda_{1}+\lambda_{2}-1)}\leq \hspace{2cm}
\\ \nonumber
(a_{11}+a_{21})^{\lambda_{1}}(a_{12}+a_{22})^{\lambda_{2}},a_{ij}>0,\lambda_{j}>0.
\eea
If $\lambda_{1}$ and $\lambda_{2}$ take $p$ and $\frac{1}{p}$ respectively, no matter how to evaluate $a_{ij}$, Theorem \ref{2.theo} cannot be derived from above inequality. If we forcibly make the following relation hold: $a_{11}^{\lambda_{1}}a_{12}^{\lambda_{2}}=x^{p}$, $a_{21}^{\lambda_{1}}a_{22}^{\lambda_{2}}=(1-x)^{\frac{1}{p}}$, $a_{11}+a_{21}=1$ and $a_{12}+a_{22}=1$, it is too difficult to work out $a_{ij}$, $\lambda_{i}$ because the equation set includes transcendental functions. So using the generalized H\"{o}lder inequality cannot derive Theorem \ref{2.theo} directly. When $\lambda_{1}$ and $\lambda_{2}$ take identical values, H\"{o}lder inequality degenerates into Cauchy-Schwarz equality, from which Theorem \ref{2.theo} cannot be inferred too.

Young inequality states that $a^{\lambda}b^{1-\lambda}\leq \lambda a+(1-\lambda)b$ for $a,b>0$ and $1>\lambda>0$, whose left side only has a single term, so from Young inequality Theorem \ref{2.theo} cannot be gotten. Minkowski inequality means that $(\sum_{k}|a_{k}+b_{k}|^{p})^{\frac{1}{p}}\leq (\sum_{k}|a_{k}|^{p})^{\frac{1}{p}}+(\sum_{k}|b_{k}|^{p})^{\frac{1}{p}}$ for $p\leq1$. It is apparent that Minkowski inequality cannot infer Theorem \ref{2.theo}. Ky Fan inequality is still a hot spot, which states $\frac{\prod_{k=1}^{n}a_{k}}{(\sum_{k=1}^{n}a_{k})^{n}}\leq \frac{\prod_{k=1}^{n}(1-a_{k})}{(\sum_{k=1}^{n}(1-a_{k}))^{n}}$ holds when $0<x\leq \frac{1}{2}$, apparently this inequality is also essentially different from Theorem \ref{2.theo}.

The reference \cite{1082} has more than 5000 inequalities, wherein two inequalities seem close to Theorem \ref{2.theo}, which are
\bea \label{eq.cp.1}
(\frac{1-x}{2})^{q}+(\frac{1+x}{2})^{q} \leq (\frac{1+x^{p}}{2})^{\frac{1}{p-1}}
\\ \nonumber
s.t. \hspace{0.3cm} x\in[0,1],p\in(1,2],\frac{1}{p}+\frac{1}{q}=1
\eea
and
\bea \label{eq.cp.2}
(\frac{1-x}{2})^{p} +(\frac{1+x}{2})^{p} \leq \frac{1}{2}(1+x^{p}),x\in[0,1],p\in(1,2].
\eea
Comparing \eqref{eq.cp.1} and \eqref{eq.cp.2} with Remark \ref{2.rmk}, we can see they are actually different, and the left side of Remark \ref{2.rmk} is larger than that of \eqref{eq.cp.1} and \eqref{eq.cp.2} due to $(\frac{1+x}{2})^{\frac{1}{p}}\geq(\frac{1+x}{2})^{p}$ for $p\geq1$. So \eqref{eq.cp.1} and \eqref{eq.cp.2} are also completely different from Theorem \ref{2.theo}.

Based on above comparisons, we think Theorem \ref{2.theo} is a new inequality which is independent upon known compared inequalities. The inequality combines two number pairs: one pair composes of $x$ and $1-x$, and the other pair $p$ and $\frac{1}{p}$. It is very interesting that we can find the summation of the first pair is $1$, while the product of the second pair is $1$. When the two pairs intertwine as shown in Theorem \ref{2.theo}, the upper bound is $1$. The two number pairs and the constant $1$ are pretty elegant, so we take it for granted that the inequality shown in Theorem \ref{2.theo} is elegant. Theorem \ref{1.theo} is a foundation for proving Theorem \ref{2.theo}, seeing $p$ as the variable. The associated discussion happens in the area corresponding to the 2nd pair. Theorem \ref{1.theo} sees $x$ as the variable, analysis is performed in the area corresponding to the 1st pair. So, in proving the proposed inequality, two areas associated with the two pairs have been considered.

In proving Theorem \ref{1.theo} and \ref{2.theo}, only derivative test method has been used. However, contrast to the traditional methods, in using the 1st derivative test to localize the extreme points, we do not explicitly solve the points because the corresponding equation are transcendental, and cannot be explicitly solved. So, we straightforwardly use the equation of the 1st derivative test to discuss whether the extreme is local maximal or minimal. On the other hand, we also use the 1st derivative as well as the boundary condition of a function to determine whether the function will get to $0$ within a special interval. These techniques, we think, can be useful to seek and prove new inequalities consisting of transcendental terms.

\section{Conclusions} \label{sec.con}
An elegant inequality $(x)^{p}+(1-x)^{\frac{1}{p}}<1$ for $p\geq1$ and $\frac{1}{2}\geq x\geq0$ has been found and proved, which combines two number pairs, $x$ and $1-x$ as well as $p$ and $\frac{1}{p}$. It is fascinating that the summation and product of the pairs are $1$. The inequality cannot be deduced from any compared known famous inequalities. It has transcendental term, so some new techniques have been built to prove the inequality. Several variants of the inequality have also been given, it is anticipated that they are useful.

\section*{Acknowledgment}
This work is supported by NSFC under grants 61860206007 and U19A2071, as well as the funding from Sichuan University under grant 2020SCUNG205
\bibliographystyle{ieeetr}
\bibliography{ref}
\end{document}